\documentclass[11pt]{article}

\usepackage{graphicx}
\usepackage{psfrag}
\usepackage{amsmath}
\usepackage{amssymb}
\usepackage{epsfig}
\usepackage{multicol}

\newtheorem{theorem}{Theorem}
\newtheorem{proposition}[theorem]{Proposition}
\newtheorem{lemma}[theorem]{Lemma}

\title{Boundary Slopes of 2-Bridge Links\linebreak Determine the Crossing Number}
\author{
Jim Hoste\\
Pitzer College\\
\\
Patrick D. Shanahan\\
Loyola Marymount University
}
\begin{document}
\maketitle
\begin{abstract} A {\it diagonal} surface in a  link exterior $M$ is a properly embedded, incompressible, boundary incompressible surface which furthermore has the same number of boundary components and same slope on each component of $\partial M$. We derive a formula for the boundary slope of a diagonal surface in the exterior of a 2-bridge link which is analogous to the formula for the boundary slope of a 2-bridge knot found by Hatcher and Thurston. Using this formula we show that the {\it diameter} of a 2-bridge link, that is, the difference between the smallest and largest finite slopes of diagonal surfaces, is equal to the crossing number. 
\end{abstract}
\section{Introduction}
Let $\mbox{cr}(K)$ denote the minimal crossing number of a knot $K$ in the 3-sphere, and let $\mbox{D}(K)$ be the diameter of the set of finite boundary slopes of the knot. It was conjectured by Ichihara that
\begin{equation}
2 \, \mbox{cr}(K) \ge D(K)
\label{ichiharaconjecture}
\end{equation}
for all knots $K$.  This conjecture has been proven for 2-bridge knots by Mattman, Maybrun, and Robinson \cite{MMR:2005} and for Montesinos knots with three or more tangles by Ichihara and Mizushima \cite{IM:2005}. Moreover, for alternating knots, the difference between the boundary slopes of the two checkerboard surfaces (in the reduced alternating diagram) is always twice the crossing number. Hence,
\begin{equation}
2 \, \mbox{cr}(K) = D(K)
\label{alternatingknots}
\end{equation}
for all alternating Montesinos knots.   Neither \cite{MMR:2005} nor \cite{IM:2005}, however, discuss possible extensions of statements (\ref{ichiharaconjecture}) or (\ref{alternatingknots}) to link exteriors. In this paper we do this by considering a restricted set of essential surfaces in the link exterior  which we call ``diagonal'' surfaces. Our main result, Theorem~\ref{bs for links}, provides a formula for the boundary slope of a diagonal surface of a 2-bridge link $L$ which is analogous to the formula given  by Hatcher and Thurston for the boundary slope of a 2-bridge knot.  As an application of this formula, we prove Theorem~\ref{diameter vs crossing number}, that
$$ \mbox{cr}(L) = D_{\Delta}(L)$$
where $L$ is a 2-bridge link and $D_{\Delta}(L)$ is the diameter of the finite slopes of diagonal surfaces.  In addition, if $L$ is a non-split, $n$ component, alternating link, and if both checkerboard surfaces are diagonal, then we show in Proposition~\ref{alternatinglowerbound} that  $\frac{2}{n} \, \mbox{cr}(L) \le D_{\Delta}(L)$.  This together with Theorem~\ref{diameter vs crossing number} suggests that
$$\frac{2}{n} \, \mbox{cr}(L) \ge D_{\Delta}(L),$$
with equality in the case of alternating links,
is a possible generalization of Ichihara's conjecture to non-split links.

The paper will proceed as follows. We begin by reexamining the beautiful relationship between boundary slopes of 2-bridge knots or links and minimal edge paths in diagrams of curve systems on the 4-punctured sphere developed by Hatcher and Thurston \cite{HT:1985} and by Floyd and Hatcher \cite{FH:1988}. In Section~2 we review the salient features of this theory and use results from our paper \cite{HS:2005} in order to derive a formula for the boundary slope of a diagonal surface. We apply this formula to prove Theorem~\ref{diameter vs crossing number} in Section~3. Finally, in Section~4, we discuss extensions of these ideas to $n$ component, non-split links.

\section{Boundary Slopes}
We begin with  some basic terminology.  An {\it essential} surface $S$ in a compact, orientable 3-manifold with boundary is a properly embedded surface which is both incompressible and boundary incompressible.  If the 3-manifold is the exterior of a link of $n$ components, then we can choose a preferred basis $\{ \mu_i, \lambda_i\}$ for each boundary torus $T^2_i$, $1 \le i \le n$.  The intersection of $S$ with $T^2_i$ is a collection of $k_i$ simple, closed, nontrivial, parallel curves which determine an isotopy class represented by $\mu_i^{p_i} \lambda_i^{q_i}$ for some co-prime integers $p_i$, $q_i$.  The {\it boundary slope} of $S$ on component $i$ is defined to be the ratio $p_i/q_i$. One can also consider the $2n$-tuple $(k_1 p_1, k_1 q_1,k_2 p_2, k_2 q_2, ... , k_n p_n, k_n q_n)$ which encodes the boundary slopes and the number of sheets of $S$ on each boundary component.  

According to Hatcher~\cite{H:1982}, a knot can have only a  finite number of boundary slopes. However for a link there may be infinitely many for each component.  Therefore, in order to define a diameter we restrict our attention to a special subset of surfaces. Define a {\it diagonal surface} to be an essential surface whose associated $2n$-tuple has the form $(kp,kq,kp,kq,...,kp,kq)$. That is, the boundary slope and number of sheets of $S$ on each component is the same.  To each diagonal surface we can assign the single slope $p/q$. If $L$ is a non-split link, then Floyd and Oertel~\cite{FO:1984} prove that there are  a finite number of branched surfaces in the link exterior that carry all essential surfaces. As pointed out in \cite{H:1982}, if $S_1$ and $S_2$ are two essential surfaces carried by the same branched surface, then the intersection number $\partial S_1 \cdot \partial S_2$ is equal to zero. If $S_1$ and $S_2$ are diagonal with $2n$-tuples  
 $(kp,kq,kp,kq,...,kp,kq)$ and $(jr,js,jr, js,...,jr,js)$, then
 $$ \partial S_1 \cdot \partial S_2=kj(pr-qs)=0.$$
 Hence, $S_1$ and $S_2$ have the same slope and so there can only be a finite number of diagonal boundary slopes for any non-split link. It follows that the diameter $D_\Delta(L)$ of any non-split link is finite.

\begin{figure}
\psfrag{a}{$1/1$}
\psfrag{b}{$3/4$}
\psfrag{c}{$2/3$}
\psfrag{d}{$1/2$}
\psfrag{e}{$1/3$}
\psfrag{f}{$1/4$}
\psfrag{g}{$0/1$}
\psfrag{h}{$-1/2$}
\psfrag{i}{$-1/1$}
\psfrag{j}{$1/0$}
\psfrag{k}{$2/1$}
\psfrag{l}{$3/2$}

\psfrag{x}{C}
\psfrag{y}{A}
\psfrag{z}{B}
\psfrag{w}{D}

\psfrag{u}{$D_1$}

    \begin{center}
    \leavevmode
    \scalebox{1.0}{\includegraphics{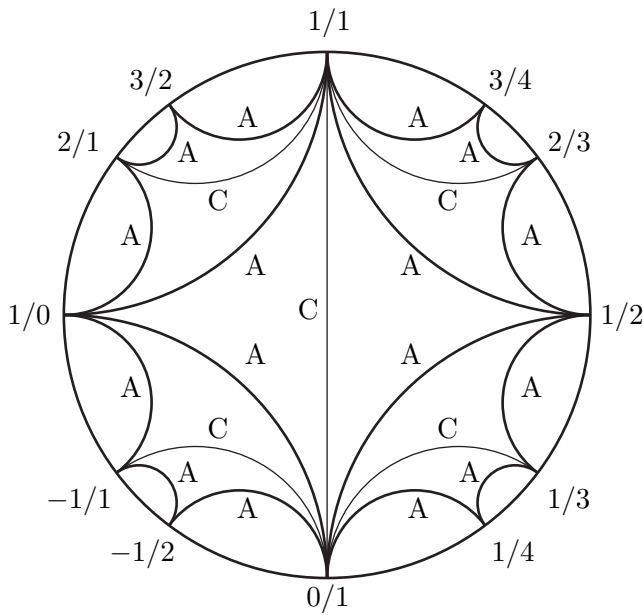}}
    \end{center}
\caption{The Diagram $D_1$.}
\label{diagram D1}
\end{figure}

We assume the reader is familiar with the fact that corresponding to each reduced rational number $p/q$ with $0\le p\le q$,  is a 2-bridge knot if $q$ is odd or 2-bridge link if $q$ is even.
However, we warn the reader that what some authors call the 2-bridge knot or link $L_{p/q}$, others call the mirror image of $L_{p/q}$. In this paper we follow the convention used in \cite{HT:1985} which is opposite that used in \cite{FH:1988}, \cite{HS:2005}, and \cite{L:1993}. In the papers by Hatcher and Thurston \cite{HT:1985} and Floyd and Hatcher \cite{FH:1988}, the set of essential surfaces in the exterior of a 2-bridge knot or link  $L_{p/q}$ are completely described and classified. These papers develop a beautiful correspondence between the essential surfaces in $L_{p/q}$ and certain paths in diagrams of curve systems on the 4-punctured sphere. For all essential surfaces in knot exteriors, and  for diagonal surfaces in link exteriors,   this diagram is called $D_1$ and is shown in Figure~\ref{diagram D1}. 

The diagram $D_1$ is a  tessellation  of $\Bbb H^2$ by ideal triangles. The rationals, together with $\frac{1}{0}$, are arranged around the unit circle as shown, and two fractions $\frac{a}{b}$ and $\frac{c}{d}$ are connected by a geodesic if and only if $ad-bc=\pm 1$.  The group of orientation preserving symmetries of $D_1$ is PSL$_2\Bbb Z$. 
Let $G\subset \mbox{PSL}_2\Bbb Z$ be the subgroup of M\"obius transformations given by $z\to \frac{az+b}{cz+d}$ where $c$ is even. It follows that the  ideal triangle $\{\frac{1}{0}, \frac{0}{1}, \frac{1}{1}\}$ is a fundamental domain for the action of $G$ and the $G$-images of the {\it fundamental} ideal quadralateral $Q=\{\frac{1}{0}, \frac{0}{1},\frac{1}{2}, \frac{1}{1}\}$ tessellate $\Bbb H^2$. There are two distinct orbits of edges which are labeled $A$ and $C$.

An edge  path in $D_1$ is {\it minimal} if it never contains two consecutive edges that lie in the same triangle. According to  \cite{HT:1985} and  \cite{FH:1988},
each minimal edge path in $D_1$ from $\frac{1}{0}$ to $\frac{p}{q}$ determines a diagonal surface in $L_{p/q}$.  (A similar correspondence exists for non-diagonal surfaces but involving paths in a more complicated diagram $D_t$.) For a particular fraction $\frac{p}{q}$ there can only be a finite number of minimal edge paths connecting it to $\frac{1}{0}$. This follows from the fact that these  minimal paths are all contained in a unique minimal chain of quadralaterals consisting of  $Q$ and a finite number of its translates under $G$. 

In order to determine the slope of a diagonal surface we must first describe several important features of edge paths in $D_1$. Each edge path from $\frac{1}{0}$ to $\frac{p}{q}$ corresponds to a continued fraction expansion\footnote{We follow the notational convention of \cite{HT:1985}.}
$$p/q=r+[b_1, b_2, \dots, b_k]=r+\frac{1}{b_1-\frac{1}
{
\begin{array}{ccc}
b_2-&&\\
&\ddots&\\
&&-\frac{1}{b_k}\\
 \end{array}
 }
 }$$
where the partial sums
$$p_i/q_i=r+[b_1, b_2, \dots, b_i]$$
are the consecutive vertices on the path. At the vertex $p_i/q_i$ the path turns left with $b_{i+1}$ triangles on the left if $b_{i+1}>0$  or to the right with $-b_{i+1}$ triangles on the right if $b_{i+1}<0$. For example, the path $\gamma$ shown in Figure~\ref{midpath with turning numbers} corresponds to the expansion
$$\frac{13}{34}=0+[2,-1,1,-1,1,-2].$$
Because $b_i$ can be interpreted in terms of the amount of turning at vertex $p_{i-1}/q_{i-1}$, we call the $b_i$'s the {\it turning numbers} of the path. For any path $\gamma$ let $n_\gamma^+$ and $n_\gamma^-$ be the number of positive and negative turning numbers respectively.  Minimality of a path can now be stated in terms of the turning numbers: a path is minimal if and only if all the turning numbers are 2 or more in absolute value.

We may recursively generate $D_1$ by starting with the initial pair $\frac{1}{0}$ and $\frac{0}{1}$ and then introducing {\it mediants}.  We first introduce the mediant $\frac{1+0}{0+1}=\frac{1}{1}$, obtaining the  sequence $\{\frac{1}{0},\frac{1}{1},\frac{0}{1} \}$. We now insert mediants again between each consecutive pair of fractions to obtain 
$\{\frac{1}{0},\frac{2}{1},\frac{1}{1},\frac{1}{2},\frac{0}{1} \}$ and so on. This process keeps the sequence in decreasing order and also preserves the fact that the determinant  $p_iq_{i+1}-p_{i+1}q_i$ of consecutive fractions is always $+1$. Viewed this way we see that every vertex in $D_1$ has two {\it parents}: the fractions that gave birth to it when taking mediants. Both parents of a link (a fraction with even denominator, and thus corresponding to a 2-bridge link) are knots, while the parents of a knot (a fraction with odd denominator) are a mixture of a knot and a link. Furthermore, the numerators of the parents of a link must have opposite parity. It is not hard to see that each vertex in a minimal path from  $\frac{1}{0}$ to $\frac{p}{q}$ must be a parent of the next vertex. 

A minimal path is called {\it even} if all of the turning numbers are even. Note that an even path starting at $1/0$ ,can never traverse the diagonal ($C$-type edge) of any quadralateral. Thus, the vertices along an even path must alternate between knots and links. Inducting on the length of a path, it is not hard to show that each knot $p/q$ has a unique even path, $e(p/q)$, connecting it to $1/0$.  Furthermore, each link has  exactly two such even paths and we denote the even path which arrives via the parent with even numerator $e^0(p/q)$ and the one which arrives via the parent with odd numerator $e^1(p/q)$. The even path to a knot is the extension of exactly one of the two even paths to its link parent.

Given $p$ and $q$ define $\epsilon_i(p/q)$ as
$$\epsilon_i(p/q)=(-1)^{\lfloor ip/q \rfloor} \quad \text{for \quad $0<i<q$},$$
where $\lfloor x \rfloor$ is the greatest integer less than or equal to $x$. 
These numbers play several important roles in relation to the 2-bridge knot or link $L_{p/q}$. For example, they can be used to express the single relation in  a certain  2-generator presentation of the fundamental group of the complement. Or, in the case of a link, the sum of all the $\epsilon_i$'s where $i$ is odd is the linking number of the two components (assuming a certain orientation convention). It is convenient to introduce notation for the sum of the even and the sum of the odd $\epsilon_i$'s.  Let
$$\sigma_0(p/q)=\sum_{i=1}^{\lfloor(q-1)/2\rfloor}\epsilon_{2 i}(p/q) \quad \text{and} \quad \sigma_1(p/q)=\sum_{i=0}^{\lfloor(q-2)/2\rfloor}\epsilon_{2i+1}(p/q).
$$

If $\gamma=\{1/0, p_1/q_1, \dots, p_n/q_n\}$ is any oriented edge-path in $D_1$ from $1/0$ to $p_n/q_n$, we define $m(\gamma)$ to be
$$m(\gamma)=\sum_{i=1}^{n-1} \det \begin{pmatrix} p_{i}&p_{i+1}\\q_{i}&q_{i+1} \end{pmatrix}.$$
More generally, if $\gamma$ is any oriented edge-path in $D_1$ we define $m(\gamma)$ to be the sum of the {\it determinants} of its edges, excluding any edge containing $\frac{1}{0}$.

We may now state the following result  which gives three different formulations for the boundary slope of a 2-bridge knot corresponding to a specific minimal path in $D_1$. The first of these is Proposition~2 of \cite{HT:1985}. Later in this section we will show that the second two formulations follow from the first.

\begin{theorem}
\label{ht prop 2} If $\gamma$ is a minimal path in $D_1$ from  $\frac{1}{0}$ to the 2-bridge  knot $\frac{p}{q}$, then each of the following gives the boundary slope of the corresponding essential surface.
\begin{enumerate}
\item $2[(n_\gamma^+-n_\gamma^-)-(n_{e(p/q)}^+-n_{e(p/q)}^-)]$
\item $-2\left [m(\gamma)-m(e(p/q))\right ]$
\item $-2[m(\gamma)-2\sigma_0(p/q)]$
\end{enumerate}
\end{theorem}

It is worth noting that the first formula for the boundary slope given above can be thought of as a  slope of $2(n_\gamma^+-n_\gamma^-)$ with respect to a non-preferred longitude, which  is then rewritten in terms of the preferred longitude by subtracting the ``correction'' term $2(n^+_{e(p/q)}-n^-_{e(p/q)})$. Similarly, $m(e(p/q))$ and $2\sigma_0(p/q)$ provide the correction in the second and third formulas.

In \cite{FH:1988}, Floyd and Hatcher extend the work of Hatcher and Thurston, classifying all essential surfaces in 2-bridge link exteriors. Lash~\cite{L:1993} then developed an algorithm to compute the boundary slopes of these surfaces. The following theorem provides the analog of Theorem~\ref{ht prop 2} for diagonal surfaces in 2-bridge link exteriors. The third of these formulations is easily derived from  \cite{HS:2005} which simplifies and extends the work of \cite{L:1993}.

\begin{theorem}
\label{bs for links}If $\gamma$ is a minimal path in $D_1$ from  $\frac{1}{0}$ to the 2-bridge  link $\frac{p}{q}$, then each of the following gives the boundary slope of the corresponding diagonal surface. 
\begin{enumerate}
\item $\displaystyle (n_\gamma^+-n_\gamma^-) -\frac{1}{2}\left [(n_{e^0(p/q)}^+-n_{e^0(p/q)}^-)+(n_{e^1(p/q)}^+-n_{e^1(p/q)}^-)\right ]$
\item $\displaystyle -\left [m(\gamma)-\frac{1}{2}\left (m(e^0(p/q))+m(e^1(p/q))\right ) \right ]$
\item $\displaystyle -[m(\gamma)-\sigma_0(p/q)]$
\end{enumerate}
\end{theorem}

The similarity between Theorem~\ref{ht prop 2} and Theorem~\ref{bs for links} can be further enhanced by thinking of the ``correction'' terms as being obtained by averaging over all even paths. For example, there are two even paths to a link and so  $(m(e^0)+m(e^1))/2$ is the average value of $m(\gamma)$ averaged over all even paths $\gamma$. Since there is only one even path to a knot, $m(e)$ is again the average value of $m$  averaged over all even paths. 

Continuing to compare Theorems~\ref{ht prop 2} and \ref{bs for links}, we also see that 
it is necessary to multiply by a factor of $2$ when going from links to knots. This makes sense because when the two components of a link are ``joined'' to form a single knot, and boundary curves on each component are connected to form a single boundary curve, the numbers of longitudes and meridians comprising each of the original boundary curves  must be combined. Since diagonal surfaces have the same data on each component this combination amounts to multiplication by $2$.

The remainder of this section will be devoted to proving Theorems~\ref{ht prop 2} and \ref{bs for links}. A number of lemmas are required to relate the quantities $n^+_\gamma-n^-_\gamma, m(\gamma)$ and $\sigma_0(p/q)$ for different paths $\gamma$. 
We begin with a result describing how $m(\gamma)$ changes when $\gamma$ undergoes a simple change.

Let $\cal C$ be the minimal chain of quadrilaterals from $\frac{1}{0}$ to $\frac{p}{q}$.       Suppose $\gamma$ is any path in $\cal C$ from $\frac{1}{0}$ to $\frac{p}{q}$ and that $T$ is a triangle in $\cal C$ having exactly one edge $e$ in $\gamma$. If we remove $e$ from $\gamma $ and replace it with the other two edges of $T$ we obtain a new path in $\cal C$ from from $\frac{1}{0}$ to $\frac{p}{q}$. This move, and its inverse, we will call a {\it  triangle move}. We may further refine our definition to {\it left} and {\it right} triangle moves depending on whether $T$ lies to the left or right of the original path. 

\begin{lemma} \label{triangle move}Changing $\gamma$ by a left triangle move increases $m(\gamma)$ by 1.
\end{lemma}
\noindent{\bf Proof:}
Every quadrilateral in $\cal C$ is the image of the fundamental quadrilateral by an element 
$$g=\left ( \begin{array}{cc} a& b \\ c& d \end{array} \right )$$
where $c$ is even  and $ad-bc=1$. Thus, the vertices of the quadrilateral are (in counter-clockwise order) 
$$\left\{ \frac{a}{c}, \frac{b}{d},  \frac{a+2 b}{c +2 d}, \frac{a+b}{c+d} \right\}.$$ Furthermore, the diagonal connects $\frac{b}{d}$ to $\frac{a+b}{c+d}$. Let $T_1$ be the triangle with vertices $\{ \frac{a}{c}, \frac{b}{d}, \frac{a+b}{c+d} \}$ and $T_2$ the triangle with vertices $\left\{ \frac{b}{d}, \frac{a+2 b}{c+2 d}, \frac{a+b}{c+d} \right\}$.  If we orient the boundary of each triangle counter-clockwise, then it is a simple matter to check that $m(\partial T_i)=-1$ for $i=1,2$. Now suppose that $\gamma^\prime$ is obtained from $\gamma$ by a left triangle move across the triangle $T$.  If $m(\gamma)=x+y$ where $y$ is the contribution to $m(\gamma)$ due to the edges of $T$ in $\gamma$, then the edges of $T$ in $\gamma'$  contribute $y+1$ to $m(\gamma')$. Hence
$$m(\gamma')=x+y+1=m(\gamma)+1.$$
Therefore, a left triangle move always increases $m$ by 1.
\hfill $\square$

 There are two paths in $\cal C$ from $\frac{1}{0}$ to $\frac{p}{q}$ which we will call the {\it upper} and {\it lower} paths. Topologically, $\cal C$ is a disk. The lower path follows the perimeter of $\cal C$ from $\frac{1}{0}$ to $\frac{p}{q}$ in the counter-clockwise direction while the upper path follows the perimeter in the clockwise direction. Except when $q=1$, neither path can contain three edges in a row from a single quadrilateral since  $\cal C$ is minimal. However, it is possible that two edges in a row are from the same quadrilateral. In this case, if the vertex common to the two edges has an even denominator, then the path is not minimal. If we replace each such occurrence with the diagonal of that quadrilateral, then the path will be minimal. Call these two paths the {\it lower minimal path} $\gamma_{\ell}$ and the {\it upper minimal path} $\gamma_u$. (If $q=1$, then $p/q=0/1$, or $1/1$ and  $\gamma_{\ell}=\gamma_u$.)

\begin{lemma} \label{lower path} The determinant of every edge in both the lower path and the lower minimal path of $\cal C$ (except for the first edge, which contains $\frac{1}{0}$) is $-1$. The determinant of every edge in both  the upper path and the upper minimal path of $\cal C$ (except for the first edge, which contains $\frac{1}{0}$) is $+1$. 
\end{lemma}
\noindent{\bf Proof:}  We show first that as we traverse the perimeter of $\cal C$ in the counter-clockwise direction the determinant of each edge, other than the edge $\{\frac{1}{0}, \frac{0}{1}\}$, is $-1$. If $\cal C$ consists of  a single quadrilateral then this is easy to check. Proceeding by induction, imagine that the last quadrilateral of the chain has been attached to all the previous ones along the edge $\{\frac{a}{c}, \frac{b}{d}\}$, where the vertex $\frac{a}{c}$ is reached before $\frac{b}{d}$ as one travels counter-clockwise from $\frac{1}{0}$. Thus, by our inductive hypothesis, $ad-bc=-1$. If $c$ is even, then the perimeter of  $\cal C$ has been changed by replacing the edge $\{\frac{a}{c}, \frac{b}{d}\}$ with the sequence of three edges $\{ \{\frac{a}{c}, \frac{a+b}{c+d}\},   \{\frac{a+b}{c+d}, \frac{a+2 b}{c+2d}\},  \{\frac{a+2 b}{c+2 d}, \frac{b}{d}\}            \}$. Each of these three new edges has determinant $-1$. If instead $d$ is even, the edge $\{\frac{a}{c}, \frac{b}{d}\}$ is replaced with the sequence    $\{ \{\frac{a}{c}, \frac{2a+b}{2c+d}\},   \{\frac{2a+b}{2c+d}, \frac{a+ b}{c+d}\},  \{\frac{a+ b}{c+ d}, \frac{b}{d}\}            \}$. Once again, the determinants of these new edges are all $-1$. Since reversing the direction of an edge negates its determinant, we see that every edge of the upper path has a determinant of $+1$.

If the lower path is not minimal, then we  may change the lower path to the lower minimal path by left triangle moves where, moreover, each triangle move replaces two edges with one edge. Such a move increases $m$ by $1$ and hence the new edge still has a determinant of $-1$. A similar argument applies to the upper minimal path.
\hfill $\square$

\begin{lemma}\label{m is n minus n} Let $\gamma$ be any path  from $1/0$ to $p/q$ in  $\cal C$. Then
$$-m(\gamma)=n_\gamma^+-n_\gamma^-.$$
\end{lemma}
\noindent{\bf Proof:} Our strategy is to first prove the result for a specific path $\gamma$ and then show that it remains true as $\gamma$ is changed to any other path by triangle moves.
Let $\gamma$ be the lower (not necessarily minimal) path in  $\cal C$.  It follows that $r=0$ and all the turning numbers $b_1, b_2, \dots, b_n$ are positive. Thus $n_\gamma^+-n_\gamma^-=n-0=n$. But, from Lemma~\ref{lower path}, $m(\gamma)=-n$. Thus,   $-m(\gamma)=n_\gamma^+-n_\gamma^-$.

Now suppose that $\gamma$ is any path in $\cal C$ and we change $\gamma$ by a left triangle move. By Lemma~\ref{triangle move} this will increase $m$ by one. We wish to show that $n_\gamma^+-n_\gamma^-$ decreases by one. Consider first the case where one edge of $\gamma$ is replaced by two edges. The new path has one more vertex and one more turning number. It is not hard to see that the new turning is negative while all the other turning numbers keep the same sign. Hence,  $n^-_\gamma$ increases by one and the difference $n_\gamma^+-n_\gamma^-$ decreases by one.

If the left triangle move exchanges two edges for one, then we may treat it as a right triangle move that exchanges one edge for two. The proof is now nearly the same as before except that the new turning number contributes to $n_\gamma^+$ instead of $n_\gamma^-$. Therefore such a right triangle move increases $n_\gamma^+-n_\gamma^-$ by one while $m$ decreases by one.
\hfill $\square$

We now turn our attention to the $\epsilon_i$'s in order to relate $\sigma_0$ to $m$ and $n^+-n^-$.

\begin{figure}
\psfrag{a}{$\epsilon=1$}
\psfrag{b}{$\epsilon=-1$}
\psfrag{p}{$q$}    
\psfrag{q}{$p$}
\begin{center}
    \leavevmode
    \scalebox{1.0}{\includegraphics{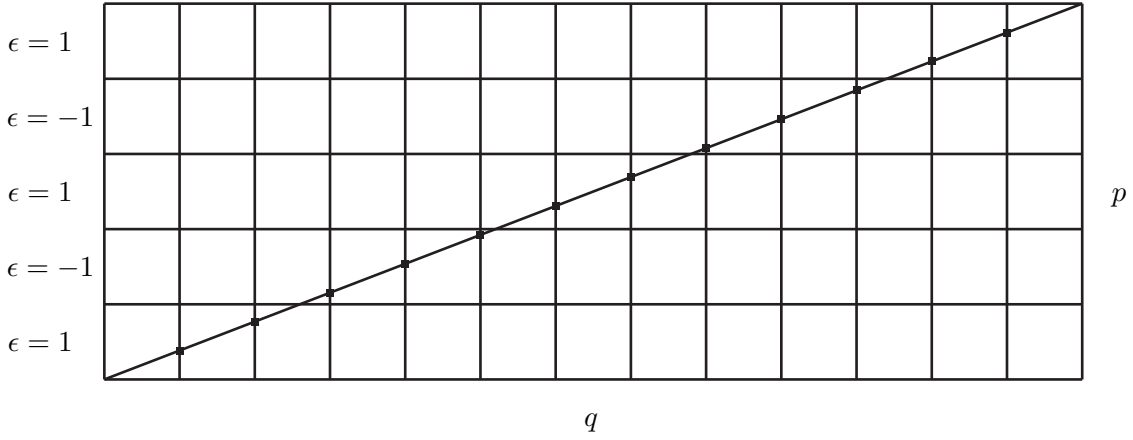}}
    \end{center}
\caption{The heights of the dots on the line of slope $p/q$ give the $\epsilon_i$'s.}
\label{epsilon graph}
\end{figure}

\begin{lemma}
\label{palindrome}
If  $0<i<q-1$, then
$$\epsilon_i(p/q)=(-1)^{p+1} \epsilon_{q-i}(p/q).$$

\end{lemma} 
\noindent{\bf Proof:} There is a beautiful, and quite useful, geometric interpretation of the $\epsilon_i$'s. Figure~\ref{epsilon graph} shows a line of slope $p/q$ extending from $(0,0)$ to $(q,p)$. It cuts the line $x=i$ at a point $P_i$ with height $ip/q$. Thus $\lfloor ip/q \rfloor$ is the height of the integer lattice point just beneath $P_i$. Each time the line passes through another horizontal line in the lattice, the signs of the $\epsilon_i$'s change. The result now follows if we consider rotating this figure 180 degrees around its center.
\hfill $\square$

The following result follows directly from Lemma~\ref{palindrome} and also can be visualized nicely in Figure~\ref{epsilon graph}.
\begin{lemma}
\label{even sum vs odd sum for knots}
If $q$ is odd, then
$$\sigma_0(p/q)=(-1)^{p+1}\sigma_1(p/q).$$
\end{lemma}

If  $p/q=(a+b)/(c+d)$ is the mediant of $a/c$ and $b/d$, then we would like to relate the sum of the even or odd $\epsilon_i$'s for $p/q$ to the corresponding sums for its parents $a/c$ and $b/d$. The next lemma provides the first step in this direction.

\begin{figure}
\psfrag{a}{$a$}
\psfrag{b}{$b$}
\psfrag{c}{$c$}
\psfrag{d}{$d$}
\psfrag{0}{$(0,0)$}
\psfrag{q}{$(q,p)$}
    \begin{center}
    \leavevmode
    \scalebox{1.0}{\includegraphics{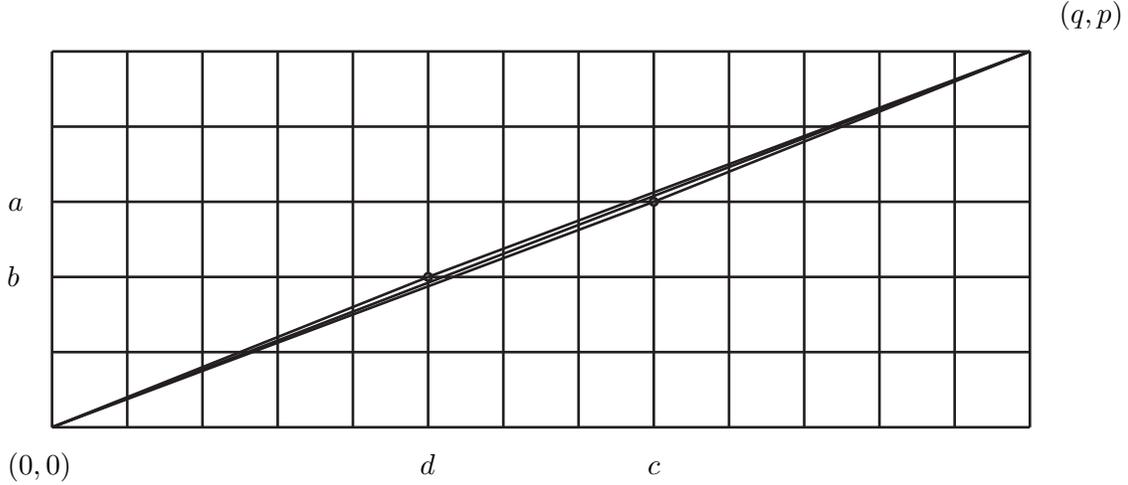}}
    \end{center}
\caption{The parallelogram contains no lattice points in its interior.}
\label{parallelogram figure}
\end{figure}

\begin{lemma}
\label{parallelogram}
Let $a,b,c$, and $d$ be positive integers such that $\gcd(a,c)=\gcd(b,d)=1$, $0<a/c<b/d$ and $ad-bc=-1$.  Let $p/q=(a+b)/(c+d)$ be the mediant of $a/c$ and $b/d$. Then
\begin{eqnarray*}
\epsilon_i(a/c)&=&\epsilon_i(p/q) \quad \mbox{ for $0<i<c$, and}\\
\epsilon_i(b/d)&=&\epsilon_i(p/q) \quad  \mbox{ for $0<i<d$.}
\end{eqnarray*}

\end{lemma}
\noindent{\bf Proof:} We give a geometric proof based on Figure~\ref{parallelogram figure}. Consider the parallelogram $P$ which is the image of the unit square, $[0,1]\times [0,1]$ under the linear transformation $T$ given by 
$$T=\begin{pmatrix}c && d\\ a && b\end{pmatrix}.$$
No integer lattice point lies in the interior of $P$ since $T$ takes the interior of the unit square to the interior of $P$.  For $0<i<c$, the points $(i,ia/c)$ and $(i,ip/q)$ lie in $P$ and hence cannot have an integer lattice point between them. Thus 
$\lfloor ia/c\rfloor=\lfloor ip/q \rfloor$ and so
$\epsilon_i(a/c)=\epsilon_i(p/q).$
A similar argument shows that 
$\epsilon_i(b/d)=\epsilon_i(p/q)$
if $0<i<d$.

(The fact that the parallelogram $P$ has an area of  1 unit is the basis of a neat parlor trick!  See page 96 of \cite{G:1961}.)
\hfill $\square$

Using Lemma~\ref{parallelogram} we may now express $\sigma_i(p/q)$ in terms of its parents. The following formulae  can all be discovered by examining Figure~\ref{parallelogram figure}.

\begin{lemma}
\label{even and odd sums of children in terms of parents}
Assume the hypotheses of Lemma~\ref{parallelogram}.

If $q$ is even, then
\begin{eqnarray*}
\sigma_0(p/q)&=&\sigma_0(a/c)+\sigma_0(b/d), \  \mbox{ and}\\
\sigma_1(p/q)&=&\sigma_1(a/c)+\sigma_1(b/d)+(-1)^a.
\end{eqnarray*} If $q$ is odd, then
$$\sigma_0(p/q)=\left\{
 \begin{array}{ll}\sigma_0(a/c)+(-1)^{p+1}\sigma_1(b/d)&\text{if $c$ is odd,}\\
 \sigma_0(b/d)+(-1)^{p+1}\sigma_1(a/c)&\text{if $d$ is odd.}
 \end{array} \right.$$
\end{lemma}
\noindent{\bf Proof:} Suppose $q$ is even and therefore both $c$ and $d$ are odd. Using Lemmas~\ref{palindrome} and \ref{parallelogram} we obtain
\begin{eqnarray*}
\sigma_0(p/q)&=&\sum_{i=1}^{(q-2)/2}\epsilon_{2i}(p/q)\\
&=&\sum_{i=1}^{(c-1)/2}\epsilon_{2i}(p/q)+\sum_{i=(c+1)/2}^{(q-2)/2}\epsilon_{2i}(p/q)\\
&=&\sum_{i=1}^{(c-1)/2}\epsilon_{2i}(a/c)+\sum_{i=(c+1)/2}^{(q-2)/2}\epsilon_{q-2i}(p/q)\\
&=&\sigma_0(a/c)+\sum_{i=1}^{(d-1)/2}\epsilon_{2i}(p/q)\\
&=&\sigma_0(a/c)+\sum_{i=1}^{(d-1)/2}\epsilon_{2i}(b/d)\\
&=&\sigma_0(a/c)+\sigma_0(b/d).
\end{eqnarray*}
If we consider the sum of the odd $\epsilon_i$'s instead, we obtain
\begin{eqnarray*}
\sigma_1(p/q)&=&\sum_{i=0}^{(q-2)/2}\epsilon_{2i+1}(p/q)\\
&=&\sum_{i=0}^{(c-3)/2}\epsilon_{2i+1}(p/q)+\epsilon_c(p/q)+\sum_{i=(c+1)/2}^{(q-2)/2}\epsilon_{2i+1}(p/q)\\
&=&\sum_{i=0}^{(c-3)/2}\epsilon_{2i+1}(a/c)+(-1)^{\lfloor cp/q \rfloor}+\sum_{i=(c+1)/2}^{(q-2)/2}\epsilon_{q-2i-1}(p/q)\\
&=&\sigma_1(a/c)+(-1)^a+\sum_{i=0}^{(d-3)/2}\epsilon_{2i+1}(p/q)\\
&=&\sigma_1(a/c)+(-1)^a+\sum_{i=0}^{(d-3)/2}\epsilon_{2i+1}(b/d)\\
&=&\sigma_1(a/c)+(-1)^a+\sigma_1(b/d).
\end{eqnarray*} The cases when $q$ is odd are similar and are left to the reader.
\hfill $\square$

\begin{proposition}\label{m vs sigma}
If $q$ is even, then
\begin{eqnarray*}
m(e^0(p/q))&=&\sigma_0(p/q)-\sigma_1(p/q), \ \mbox{and}\\
m(e^1(p/q))&=&\sigma_0(p/q)+\sigma_1(p/q).
\end{eqnarray*}
If $q$ is odd, then
$$m(e(p/q))=2\sigma_0(p/q).$$
\end{proposition}

\noindent{\bf Proof:} We induct on the number of quadrilaterals in $\cal C$. If $\cal C$ is a single quadrilateral, then the proposition is easily verified.

\begin{figure}[b]
\psfrag{p}{$P$}
\psfrag{q}{$Q$}
\psfrag{r}{$R$}
\psfrag{s}{$S$}

    \begin{center}
    \leavevmode
    \scalebox{1.0}{\includegraphics{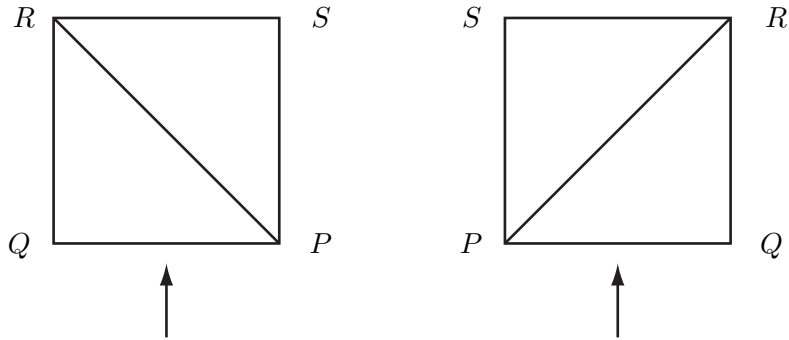}}
    \end{center}
\caption{The last quadrilateral in the minimal chain from $1/0$ to either $R$ or $S$.}
\label{last quad}
\end{figure}

To prove the inductive step, we must consider two cases depending on how the last quadrilateral is attached to the chain. The two cases are shown in Figure~\ref{last quad}. Here $P, Q, R$ and $S$ represent reduced fractions with $R$ the mediant of $P$ and $Q$, and $S$ the mediant of $P$ and $R$. The denominators of $Q$ and $S$ are even and these vertices correspond to links. The opposite is true of $P$ and $R$. Finally, the numerators of  $P$ and $R$ have opposite parity. In both cases the arrow is used to indicate that the quadrilateral is attached to the previous quadrilateral in the chain along the edge $PQ$.
 
Consider the case shown on the left side of Figure~\ref{last quad}.  We first show that the result holds for the vertex $R$ and then use this to prove the result for the vertex $S$. Throughout the proof we shall use the notation $[n(A)]$ to denote the numerator of vertex $A$ modulo 2.   From Lemma~\ref{even and odd sums of children in terms of parents} we have
$$2 \sigma_0(R) = 2 \sigma_0(P)+(-1)^{n(R)+1}2 \sigma_1(Q).$$
Using our inductive hypothesis for vertex $P$ and the fact the $n(R)$ and $n(P)$ have opposite parity we can rewrite this equation as
\begin{equation}
\label{eq A}
2 \sigma_0(R) = m(e(P))+(-1)^{n(P)}2 \sigma_1(Q).
\end{equation}
The unique even path $e(P)$ to $P$ is the extension of an even path to one of its parents and there are two cases to consider depending on whether  $P$ is a parent of $Q$ or vice versa. Suppose first  that $P$ is a parent of $Q$. The two even paths to $Q$ arrive through its parents. Thus by definition we have that $e^{[n(P)]}(Q)$ is obtained from $e(P)$ by adjoining the edge $PQ$. Since the determinant of edge $PQ$ is $-1$ we have
\begin{equation}
\label{eq B}
m(e^{[n(P)]}(Q))=m(e(P))-1.
\end{equation}
There are two even paths to $Q$ and the even path to $R$ is the extension of one of them. It cannot be the extension of $e^{[n(P)]}(Q)$ since in this case the last turning number would be $-1$. Hence it must be the extension of $e^{[n(R)]}(Q)$. (Here we use the fact that $P$ and $R$ have numerators of opposite parity). Thus
\begin{equation}
\label{eq C}
m(e(R))=m(e^{[n(R)]}(Q))+1
\end{equation}
since the determinant of $QR$ is $+1$. If instead, $Q$ is the parent of $P$ then it turns out that Equations~\ref{eq B} and \ref{eq C} are still true, an exercise that we leave to the reader.

Thus, in any case, we see that $m(e(P)) = m(e^{[n(P)]}(Q))+1$. Substituting this into Equation~\ref{eq A}, using the inductive hypothesis for vertex $Q$ twice, and finally using Equation~\ref{eq C} we obtain:
\begin{eqnarray*}
2 \sigma_0(R) & = & m(e^{[n(P)]}(Q))+1+(-1)^{n(P)}2 \sigma_1(Q) \nonumber \\ 
& = & \sigma_0(Q)-(-1)^{n(P)}\sigma_1(Q)+1+(-1)^{n(P)}2 \sigma_1(Q)  \nonumber \\
& = & \sigma_0(Q)+(-1)^{n(P)}\sigma_1(Q)+1 \nonumber \\
& = & m(e^{[n(P)+1]}(Q))+1\\
&=&m(e(R)).
\end{eqnarray*}
Hence the result is true for vertex $R$. Now consider vertex $S$ in the left-hand side of Figure~\ref{last quad}. From Lemma~\ref{even and odd sums of children in terms of parents} we have:
\begin{eqnarray*}
\sigma_0(S) & = & \sigma_0(P) + \sigma_0(R) \quad \text{and} \\
\sigma_1(S) & = & \sigma_1(P) + \sigma_1(R) + (-1)^{n(P)}.
\end{eqnarray*}
Adding these two equations and using Lemma~\ref{even sum vs odd sum for knots}, the inductive hypothesis, and the result we have already proven for vertex $R$,  we obtain
\begin{eqnarray*}
\sigma_0(S)+\sigma_1(S)&=&\left \{ \begin{array}{ll} 2 \sigma_0(R)+1& \text{if $n(P)$ is even;}\\
2 \sigma_0(P)-1&\text{if $n(P)$ is odd.} \end{array}
\right. \\
&=&\left \{ \begin{array}{ll} m(e(R))+1& \text{if $n(P)$ is even;}\\
m(e(P))-1&\text{if $n(P)$ is odd.} \end{array}
\right. \\
&=&m(e^1(S)),
\end{eqnarray*}
since the determinant of $RS$ is $+1$ and the determinant of $PS$ is $-1$.

If instead we subtract, we obtain
\begin{eqnarray*}
\sigma_0(S)-\sigma_1(S)&=&\left \{ \begin{array}{ll} 2 \sigma_0(P)-1& \text{if $n(P)$ is even;}\\
2 \sigma_0(R)+1&\text{if $n(P)$ is odd.} \end{array}
\right. \\
&=&\left \{ \begin{array}{ll} m(e(P))-1& \text{if $n(P)$ is even;}\\
m(e(R))+1&\text{if $n(P)$ is odd.} \end{array}
\right. \\
&=&m(e^0(S)).
\end{eqnarray*}
The case for the right-hand side of Figure~\ref{last quad} is similar.
\hfill $\square$

 We may now prove Theorems \ref{ht prop 2} and \ref{bs for links}.

{\bf Proof of Theorem \ref{ht prop 2}:}  Part 1 of Theorem~\ref{ht prop 2} is Proposition~2 of \cite{HT:1985}.  Part~2 now follows from Lemma~\ref{m is n minus n} and part~3 follows from Proposition~\ref{m vs sigma}. 
\hfill $\square$

{\bf Proof of Theorem \ref{bs for links}:}  Part~3 of Theorem \ref{bs for links} may be derived from \cite{HS:2005} as follows.  Lemma~3 of \cite{HS:2005} states that the boundary slope of a diagonal surface corresponding to a minimal edge-path $\gamma$ with no $C$-type edges is $m(\gamma)$.  To express this slope with respect to a preferred longitude we must subtract $\sigma_0(p/q)$ as described at the end of Section~3 in \cite{HS:2005}.  Finally, as mentioned already, the definition of $L_{p/q}$ used in \cite{HS:2005} is the mirror image of what is used here.  Thus, the boundary slope is $- \left[ m(\gamma)-\sigma_0(p/q) \right]$.

If $\gamma$ contains $C$-type edges, then Theorem~6 of \cite{HS:2005} applies.  For diagonal surfaces, the theorem gives a boundary slope of $x-P+N$.  Here $x = m(\gamma')$ where $\gamma'$ is a path with no $C$-type edges obtained from $\gamma$ by $P$ left triangle moves and $N$ right triangle moves. Hence, by Lemma~\ref{triangle move}, $x-P+N = m(\gamma)$.  Again, we must subtract $\sigma_0(p/q)$ and negate the result.

Part~2 of Theorem~\ref{bs for links} now follows from Proposition~\ref{m vs sigma}.  Finally, part~1 follows from Lemma~\ref{m is n minus n}.
\hfill $\square$
\section{Diameter and Crossing Number}
If  $L$ is a link,  let 
$D_{\Delta}(L)$ be the {\it diameter} given by the difference between the maximum and minimum (finite) slopes of diagonal surfaces in $L$. If in fact $L$ is a knot this reduces to the usual notion of diameter. 
Finally,  let $\mbox{cr}(L)$ denote the crossing number of $L$. The results of the previous section now allow us to relate the diameter of either a 2-bridge knot or link to its crossing number. We do this in the following theorem, which in the case of knots was proven in \cite{MMR:2005} using different techniques.

\begin{theorem} \label{diameter vs crossing number}If $L$ is a 2-bridge knot or link with $n$ components, then $$D_{\Delta}(L)=\frac{2}{n}\,{\rm cr}(L).$$
\end{theorem}

\noindent{\bf Proof:}  Suppose $L=L_{p/q}$ is any 2-bridge knot or link and $\cal C$ is the minimal chain of quadrilaterals from $\frac{1}{0}$ to $\frac{p}{q}$.   Let  $\gamma_\ell$ and $\gamma_u$ be the lower and upper minimal paths in $\cal C$ respectively. From Lemma~\ref{triangle move} we see that $m(\gamma_\ell)$ and $m(\gamma_u)$ provide the extreme values of $m$ since left triangle moves increase $m$ by one. From Theorems~\ref{ht prop 2} and \ref{bs for links} it follows that the diameter is
$$
D_\Delta(L)=\frac{2}{n}\left [m(\gamma_u)-m(\gamma_\ell)\right ].
$$
 It remains to show that $m(\gamma_u)-m(\gamma_\ell)$ equals the crossing number.

\begin{figure}
\psfrag{a}{$\frac{1}{0}$}
\psfrag{b}{$\frac{0}{1}$}
\psfrag{c}{$\frac{1}{1}$}
\psfrag{d}{$\frac{1}{2}$}
\psfrag{e}{$\frac{1}{3}$}
\psfrag{f}{$\frac{2}{5}$}
\psfrag{g}{$\frac{3}{8}$}
\psfrag{h}{$\frac{8}{21}$}
\psfrag{i}{$\frac{5}{13}$}
\psfrag{j}{$\frac{13}{34}$}
\psfrag{u}{$\cal U$}
\psfrag{l}{$\cal L$}
\psfrag{g1}{$\gamma_u$}
\psfrag{g2}{$\gamma_\ell$}
\psfrag{g3}{$\gamma$}

    \begin{center}
    \leavevmode
    \scalebox{.70}{\includegraphics{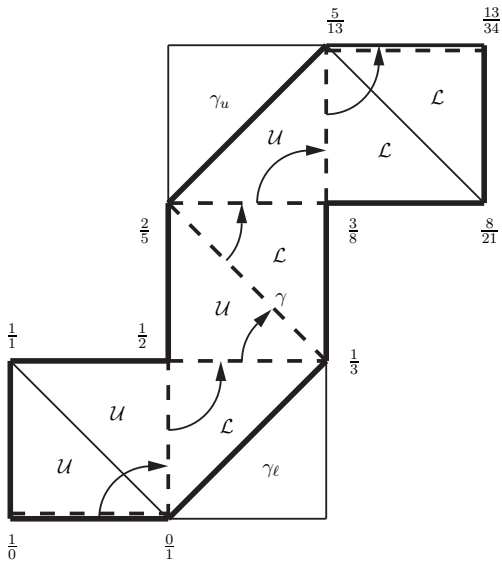}}
    \end{center}
\caption{The path $\gamma$ has turning numbers $2,-1,1,-1,1,-2$.}
\label{midpath with turning numbers}
\end{figure}

We illustrate the idea of the proof of this fact with the example $\frac{p}{q}=\frac{13}{34}$ shown in Figure~\ref{midpath with turning numbers}.

From Lemma~\ref{lower path} we have that $m(\gamma_u)=i-1$ and $m(\gamma_\ell)=-j+1$ where $i$ and $j$ are the number of edges respectively in the two paths. Thus $m(\gamma_u)-m(\gamma_\ell)=i+j-2$.

The area between $\gamma_u$ and $\gamma_\ell$ is made up of triangles which, except for the first and last triangle, have one edge on one path and the opposite vertex on the other path. Label the first triangle $\cal U$, the last triangle $\cal L$ and all the intermediate triangles either $\cal U$ or $\cal L$ depending on whether or not they contain an edge on the upper or lower path. This labeling determines a unique path $\gamma$ which keeps all the triangles labeled $\cal U$ on its left and all the triangles labeled $\cal L$ on its right. The path $\gamma$ is shown dashed in Figure~\ref{midpath with turning numbers}. Note that, as in this example,  $\gamma$ may not be minimal. The turning numbers for $\gamma$ alternate in sign and the sum of the absolute values of these turning numbers  is clearly equal to the number of triangles that have been labeled. Since each triangle contains one edge of $\gamma_u$ or $\gamma_\ell$, except for the first and last,  this number equals $i+j-2$. 
Finally,  it is well known that $\gamma$ corresponds to the continued fraction which gives the standard alternating 4-plat diagram of the link. Thus, the sum of the absolute values of the turning numbers is exactly the crossing number of $L$.
\hfill $\square$  
\section{Final Remarks}
In the introduction it was mentioned that the difference between the slopes of the two checkerboard surfaces in the reduced alternating diagram of an alternating knot is equal to twice the crossing number. That is, for an alternating knot $K$ we have $D_\Delta(K) \ge 2\, \mbox{cr} (K)$.  However, this inequality does not immediately  generalize to one for alternating 2-component links.  While \cite{DR:1999} guarantees that checkerboard surfaces in reduced, alternating diagrams of non-split  links are essential, they may not be diagonal.  For example, one of the checkerboard surfaces for the Whitehead link has slopes of $-4$ and $-2$ on the two components.  On the other hand, there are infinitely many examples of alternating links for which both checkerboard surfaces are diagonal. (For example, the three-component pretzel link (3,2,3,2,3,2) has diagonal checkerboard surfaces with slopes $-2$ and $8$, and the three-component Montesinos link $K(1/2,13/17,1/2,3/5,1/10,3/5)$ has diagonal checkerboard surfaces with slopes $-10$ and $10$).  For alternating links where both checkerboard surfaces are diagonal we have the following result.

\begin{figure}[h]
\psfrag{a}{ $K_j$}
\psfrag{b}{$K_i$}
\begin{center}
\begin{center}

    \leavevmode

    \scalebox{1.0}{\includegraphics{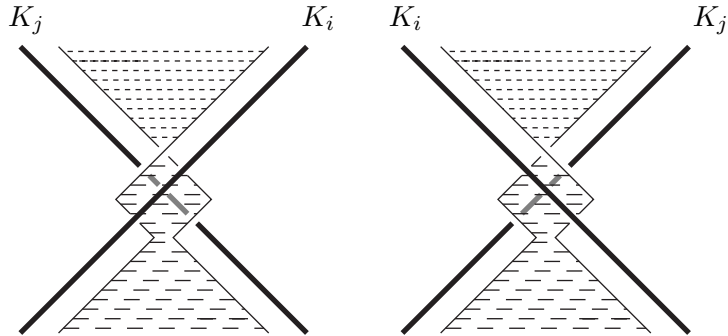}}

    \end{center}
\caption{Right and left twisted bands in a checkerboard surface.}
\label{bands}
\end{center}
\end{figure}

\begin{proposition} Let $L$ be a non-split alternating link of $n$ components, and assume that both checkerboard surfaces in a reduced alternating diagram of $L$ are diagonal.  Then
$$D_{\Delta}(L) \ge \frac{2}{n} \, {\rm cr}(L).$$
\label{alternatinglowerbound}
\end{proposition}
The proof of this proposition makes use of the following lemma.
\begin{lemma} 
\label{twisted band lemma}Let $L$ be a non-split alternating link of $n$ components.  Let $S$ and $T$ be the two checkerboard surfaces in a reduced alternating diagram of $L$.  If $s_i$ and $t_i$ are the boundary slopes of $S$ and $T$ respectively on the the $i$-th component, then
$$\left| \sum_{i=1}^n (s_i-t_i) \right| = 2 \, \mbox{\rm cr}(L).$$
\end{lemma}
{\bf Proof:}  By \cite{DR:1999} both $S$ and $T$ are essential surfaces. The surface $S$ is a collection of non-nested planar disks connected to each other by twisted bands as shown in Figure~\ref{bands}.
By examining any such disk of $S$, it is easy to see that all of the twisted bands are the same handedness because $L$ is alternating. Furthermore, since the disks of $T$ are the complementary planar regions of the disks of $S$, it follows that the bands of $T$ have the opposite handedness.  Without loss of generality, assume that $S$ has only right twisted bands while $T$ has only left twisted bands. From Figure~\ref{bands} the contribution to the boundary slope on component $i$ is easily computed.  The following table summarizes these contributions.

$$\begin{tabular}{|l|l|c|} \hline
band type & crossing type & contribution \\ \hline
 & $i \neq j$ & $+1$ \\ \cline{2-3}
right twisted & $i=j$, positive & $+2$ \\ \cline{2-3}
 & $i=j$, negative & $\ \  \, 0$ \\ \hline
 & $i \neq j$ & $-1$ \\ \cline{2-3}
left twisted & $i=j$, positive & $\ \ \, 0$ \\ \cline{2-3}
 & $i=j$, negative & $-2$ \\ \hline
\end{tabular}$$

For $1 \leq i \leq n$, let $\alpha_i$ be the number of crossings in the reduced diagram where component $i$ passes over a different component. For self-crossings,  let $P_i$ and $N_i$ be the number of positive and negative self-crossings respectively for component $i$.  Using the table we see that $s_i = \alpha_i + 2 P_i$ and $t_i = -\alpha_i - 2 N_i$.  Therefore,
$$\sum_{i=1}^n (s_i-t_i)  = 2\sum_{i=1}^n (\alpha_i+P_i+N_i) = 2 \, \mbox{cr}(L).$$
\hfill $\Box$

The proof of Proposition~\ref{alternatinglowerbound} is now simple. Since both $S$ and $T$ are diagonal, $s_i=s_j=s$ and $t_i = t_j =t$ for all $1 \le i, j \le n$.  Therefore, by Lemma~\ref{twisted band lemma} we have $n (s-t) = 2 \, \mbox{cr}(L)$ and the result follows.

\bibliographystyle{../hplain}
\bibliography{../hs}

 \end{document}